\documentclass[12pt,reqno]{amsart} 
\usepackage{amsthm}
\usepackage{amsmath,amssymb}
\usepackage{ascmac}
\usepackage{mathrsfs}
\usepackage[abbrev]{amsrefs}
\usepackage{graphicx,xcolor}
\newcommand{\BLACK}{\color{black}}

\usepackage[%
 dvipdfmx,
 setpagesize=false,%
 bookmarks=true,%
 bookmarksdepth=tocdepth,%
 bookmarksnumbered=true,%
 colorlinks=false,%
 pdftitle={},%
 pdfsubject={},%
 pdfauthor={},%
 pdfkeywords={}%
]{hyperref}            
\allowdisplaybreaks[1]
\theoremstyle{definition}
\newtheorem{dfn}{Definition}
\newtheorem{thm}{Theorem}[section]

\newtheorem{prop}[thm]{Proposition}
\newtheorem{lem}[thm]{Lemma}

\newtheorem*{nota}{Notation}

\newcommand{\pt}{\partial}             
\renewcommand{\th}{\theta}

\newcommand{\F}{\mathcal {F}}          
\renewcommand{\S}{\mathcal {S}}

\newcommand{\Z}{\mathcal {Z}}

\newcommand{\re}{\mathbb R}

\newcommand{\Int}{\mathbb Z}

\renewcommand{\L}{\Delta}
\newcommand{\Lm}{\Lambda}

\newcommand{\del}{\delta}

\newcommand{\x}{\xi}

\renewcommand{\t}{\tau}

\newcommand{\N}{\nabla }

\renewcommand{\div}{{\rm {div\ }}}

\def\<{\langle }

\newcommand{\supp}{\text{ supp }}

\newcommand{\eqsp}[1]{{\begin{equation}\begin{split}#1\end{split}\end{
equation}}}
\begin{document}

\begin{center}\LARGE \bf
 On the uniqueness of the mild solution of the critical quasi-geostrophic equation
\end{center}
   
\footnote[0]
{


E-mail: *t-iwabuchi@tohoku.ac.jp, **okazaki.taiki.r5@dc.tohoku.ac.jp

}
\vskip5mm

\begin{center}
 {\large Tsukasa Iwabuchi* and Taiki Okazaki**} 
 \vskip2mm
 {\large Mathematical Institute, Tohoku University}\\
 {\large Sendai 980-8578 Japan}
\end{center}

\vskip5mm

\begin{center}
\begin{minipage}{135mm}
\footnotesize
{\sc Abstract.}
We demonstrate that the uniqueness of the mild solution of the two-dimensional  quasi-geostrophic equation with the critical dissipation
holds in the scaling critical homogeneous Besov space $\dot{B}^0_{\infty,1}$.
We consider a solustion of integral equation, and our result does not need regularity assumption.
\end{minipage}
\end{center}

\section{Introduction.}
In this paper, we consider the critical quasi-geostrophic equation in $\re^2$.

\begin{equation}
   \label{SQG}
   \begin{cases}
     \pt_t \th + \Lm \th + u \cdot \N \th = 0,
           &\quad  \   t>0,x\in{\re}^2, \\
     u = \N ^\perp \Lm^{-1}\th,
           &\quad  \ t>0,x\in{\re}^2,\\
     \th(0,x) = \th_0(x),
           &\quad  \ x\in{\re}^2,\\
   \end{cases}
\end{equation}
where $\N ^\perp = (-\pt_{x_2},\pt_{x_1})$ and $\ \Lm=(-\L)^\frac{1}{2}.$
The real-valued function $\th(t,x)$ denotes the potential temperature of the fluid.
The quasi-geostrophic equation is an important model in geophysical dynamics, 
which describes large scale atmospheric and oceantic motion with small Rossby and Ekman numbers (see \cite{Pe_1979}).
From a mathematical point of view, the quasi-geostrophic equation has similar stracture with the there-dimensional Euler and Navier-Stokes equations (see \cite{Co_1994}).

Let us recall several known results about the existence and the uniqueness. 
We consider $\Lm^\alpha$ with $1\leq \alpha \leq 2$ instead of $\Lm$ in the equation of \eqref{SQG}.
In the cace when $1< \alpha \leq 2$, Constantin and Wu \cite{Co_1999} proved that the strong solution 
$\th\in L^\infty(0,T;L^2)\cap L^2(0,T;H^\frac{\alpha}{2})\cap L^q(0,T;L^p)$ exists uniquelly for 
$\th_0\in L^2,p\geq 1,q>1,1/p+\alpha/2q=\alpha/2-1/2$.
If $1< \alpha < 2$, Ferreira \cite{Fe_2011} established the uniqueness of the mild solution $\th\in C([0,T];L^\frac{2}{\alpha-1})$.

For $\alpha=2$, Iwabuchi, Ueda \cite{I_U_2024} showed the uniqueness of mild solution $\th\in C([0,T]; L^2)$. 
If $0<\alpha\leq1$, the existence of the local unique solution $\th\in C([0,T];\mathring{B}^{1-\alpha}_{\infty,q})\cap \widetilde{L}^1_TB^1_{\infty,q}$ for $\th_0\in\mathring{B}^{1-\alpha}_{\infty,q},1\leq q<\infty$, is proved by Wang and Zhang (\cite{Wa_Zh_2011}). 
Here $\mathring{B}^{1-\alpha}_{\infty,q}$ is the closure of the Schwartz functions in the norm of the inhomogenous Besov space $B^{1-\alpha}_{\infty,q}$.
This paper is devoted to the uniqueness of mild solutions of \eqref{SQG} in a critical space $\dot{B}^0_{\infty,1}$.

We also mention some literature on the Navier-Stokes equations. 
The scale critical spaces of the two-dimensional surface quasi-geostrophic equation are same as those of the two-dimensional incompressible Navier-Stokes equations. 
It is a classical result that the weak solustion $u\in L^\infty(0,T;L^2)\cap L^2(0,T;\dot{H}^1)$ is unique if $u_0\in L^2$ for
the two-dimensional incompressible Navier-Stokes equations.
The non-uniqueness in $C(0,T;L^p)$ with $1\leq p<2$ for the two-dimensional incompressible Navier-Stokes equations
on $\mathbb{T}^2$ is provided by Cheskidov-Luo \cite{Ch_2023}.
In the case when  the space dimension is three, 
Meyer (\cite{Me_1996}), Furioli, Lemari\'{e}-Rieusset, Terraneo (\cite{Fu_1997}) and Monniaux (\cite{Mo_1999}) proved the uniqueness of the mild solution in $C([0,T];L^3)$
for the incompressible Navier-Stokes equations.
The corresponding result for our problem is the uniqueness for $\th\in C([0,T];L^\infty)$, 
but there are some difficulty. 
Hence we consider \eqref{SQG} in $\dot{B}^0_{\infty,1}$ which is close to $L^\infty$. 

Let us recall the definition of homogeneous Besov space. 
We refer to the book by Bahouri, Chemin, Danchin~(\cite{Danchin}). 
\begin{dfn}
   Let $\{\phi_j\}_{j\in\Z}\subset\S(\re^2)$\ be such that
      $$\supp \widehat{\phi_j}\subset\{\xi\in\re^2\ |\ 2^{j-1}\leq |\xi| \leq 2^{j+1}\}\  \text{for any}\ j\in\Int, \quad \sum_{j\in\Int}\widehat{\phi_j}(\xi)=1\ \text{for any}\ \xi\in\re^2\backslash\{0\}.$$
   Let $\S_{h}'$ be the space of tempred distribution $f$ such that
   $$\lim_{\lambda\to \infty}{\left\|\F^{-1}\left[\psi(\lambda\x)\widehat{f}(\x)\right]\right\|}_{L^\infty}=0\quad \text{for any}\ \psi\in C^\infty_c(\re^2).$$
   \BLACK
   Then
   $$f=\sum_{j\in\Int}\phi_j*f$$
   holds in $\S_{h}'$.
   For $s\in\re$\ and $1\leq p,q \leq \infty$,\ we define the homogeneous Besov spaces as follows.
   $$\dot{B}^s_{p,q}=\dot{B}^s_{p,q}(\re^2):=\{f\in \S_{h}'\ |\ {\|f\|}_{\dot{B}^s_{p,q}}<\infty\},$$
   where
   $${\|f\|}_{\dot{B}^s_{p,q}}:={\left\|\left\{2^{sj}{\|\phi_j*f\|}_{L^p}\right\}_{j\in\Int}\right\|}_{l^q}.$$
\end{dfn}
\noindent
{\bf Remark. } It is known that if $s<2/p$ or $s=2/p$ and $q=1$, $\dot{B}^s_{p,q}$ is a Banach space.
\vskip3mm
We define the mild solution of \eqref{SQG}.
\begin{dfn}
   Let $T>0,\th_0\in \dot{B}^0_{\infty,1}$.\ If $\th:[0,T]\times\re^2\to\re$\ satisfies
   \begin{equation}\notag
      \begin{cases}
       \th\in C([0,T];\dot{B}^0_{\infty,1}), \\
       \displaystyle\int_{\re^2}\th(t,x)\phi(x) ~{\rm d}x=\int_{\re^2}e^{-t\Lm}\th_0(x)\,\phi(x)~{\rm d}x+\int_{\re^2}\left(\int_{0}^{t}e^{-(t-s)\Lm}(u\th) ~{\rm d}s\right)\cdot\N\phi(x) ~{\rm d}x, 
      \end{cases}
   \end{equation}
   for all \ $t\in[0,T]\ \text{and} \ \phi\in\S(\re^2)$, 
   then we call $\th$ a mild solution of (\ref{SQG}).
\end{dfn}

We state our main result. 
\begin{thm}
   \label{Thm}
   Let $T>0$\ and $\theta^{(1)},\theta^{(2)}\in C([0,T];\dot{B}^0_{\infty,1})$\ are mild solutions of (\ref{SQG}) such that $\theta^{(1)}(0)=\theta^{(2)}(0)$.
   Then $\theta^{(1)}(t)=\theta^{(2)}(t)$\ in $\dot{B}^0_{\infty,1}$\ for all $t\in[0,T]$.
\end{thm}

In the proof of Theorem \ref{Thm}, we divide the equation \eqref{SQG} into high frequency part and low frequency part
 and consider the norms $\dot{B}^{-\eta}_{\infty,\infty}$ and $\dot{B}^0_{\infty,1}$ with $0<\eta<1$, respectively. \BLACK
After that using the maximum principle in homogeneous Besov space that is proved by \cite{Miao}.
To this end, we introduce some notations.

\begin{nota}
   For any $j\in\Int$, we denote
   \begin{equation}
      \label{0321-1}
      f_j:=\phi_j*f,\quad f_{\text{high}}:=\sum_{j=1}^\infty\phi_{j}*f\quad \text{and}\quad  f_{\text{low}}:=f-f_{\text{high}}.\BLACK 
   \end{equation}
   Also we introduce a notation for nonlinear term in the integral equation of (\ref{SQG}),
   \begin{equation}
      \label{0313-1}
      \th(t)=e^{-t\Lm}\th_0-\int_{0}^{t}\N\cdot e^{-(t-s)\Lm}(u\th) ~{\rm d}s=:e^{-t\Lm}\th_0-B(u\th)(t).
   \end{equation}
\end{nota}

\section{Preliminarys.}
In this section, we introduce several lemmas which are elemental properties of Besov spaces and heat kernel.

\begin{lem}(\cite{Danchin,Wu_2008})
   \label{0314-1}
   Let $s\in\re$ and $1\leq p\leq \infty$.
   Then there exist constants $C,c>0$ such that
   \begin{equation}
      \label{0314-2}
      {\|\N(\phi_j* f)\|}_{L^p}\leq C2^j{\|\phi_j*f\|}_{L^p},\ {\|\Lm^{-1}(\phi_j* f)\|}_{L^p}\leq C2^{-j}{\|\phi_j*f\|}_{L^p}
   \end{equation}
   and
   \begin{equation}
      \label{0407-1}
      {\|\phi_j*e^{-t\Lm} f\|}_{L^p}\leq Ce^{-ct2^j}{\|\phi_j*f\|}_{L^p}\ 
   \end{equation}
   for all $f\in L^p,t>0,j\in\Int$.
\end{lem}

\begin{lem}(\cite{Danchin})
   \label{0314-4}
   Let $s\in\re,\ s'<0$ , $1\leq p,q\leq\infty$ and $1/q= 1/q_1+1/q_2\leq 1$.
   Then
   $${\left\|\sum_{k\leq l-2}f_kg_l\right\|}_{\dot{B}^s_{p,q}}\leq C{\|f\|}_{L^\infty}{\|g\|}_{\dot{B}^s_{p,q}}\ 
   \text{for all}\ f\in L^\infty,\ g\in\dot{B}^s_{p,q},$$
   and
   $${\left\|\sum_{k\leq l-2}f_kg_l\right\|}_{\dot{B}^{s+s'}_{p,q}}\leq C{\|f\|}_{\dot{B}^{s'}_{\infty,q_1}}{\|g\|}_{\dot{B}^s_{p,q_2}}\ 
   \text{for all}\  f\in \dot{B}^{s'}_{\infty,q_1},\ g\in\dot{B}^s_{p,q}.$$
\end{lem}

\begin{lem}(\cite{Danchin})
   \label{0314-5}
   Let $s_1,s_2\in\re,\ s_1+s_2>0, 1\leq p,q\leq\infty, 1/p=1/p_1+1/_2\leq 1$ and $1/q=1/q_1+1/q_2\leq 1.$
   Then
   $${\left\|\sum_{|k-l|\leq 1}f_kg_l\right\|}_{\dot{B}^{s_1+s_2}_{p,q}}\leq C{\|f\|}_{\dot{B}^{s_1}_{p_1,q_1}}{\|g\|}_{\dot{B}^{s_2}_{p_2,q_2}}\ 
   \text{for all}\  f\in\dot{B}^{s_1}_{p_1,q_1},\ g\in\dot{B}^{s_2}_{p_2,q_2}.$$
\end{lem}

\begin{lem}(\cite{Danchin})
   \label{0323-1}
   Let $s>0,\ 1\leq p,q\leq\infty .$
   Then
   $${\|fg\|}_{\dot{B}^s_{p,q}}\leq C{\|f\|}_{L^\infty\cap\dot{B}^s_{p,q}}{\|g\|}_{L^\infty\cap\dot{B}^s_{p,q}}\ 
   \text{for all}\  f,g\in L^\infty\cap\dot{B}^s_{p,q}.$$
\end{lem}

Next lemma is concerned with the boundness of biliner Fourier multiplier in $L^\infty$. 
We notice the frequency of the function is restricted. 
\BLACK

\begin{lem}
   \label{0314-12}
   Let $j,j'\in\Int$ and $|j-j'|\leq 1$.
   Then
   \begin{equation}
      \label{0420-2}
      {\left\|\int_{\re^2\times\re^2}\frac{(\xi,-\eta)}{|\x|+|\eta|}\widehat{f_j}(\x)(\eta)\widehat{g_{j'}}(\eta) e^{i(\cdot,\cdot)\cdot(\xi,\eta)}\BLACK  ~{\rm d}(\x,\eta)\right\|}_{L^\infty}\leq C{\|f_j\|}_{L^\infty}{\|g_{j'}\|}_{L^\infty}
   \end{equation}
   for all $f,g\in L^\infty$, where $f_j, g_{j'}$ are defined by \eqref{0321-1}
\end{lem}
\begin{proof}
We can write left side of \eqref{0420-2} using Fourier multiplier in $\re^4$, 
which we write by $\F^{-1}_{\re^4}$. 
We denote $\widetilde{\phi_j}=\phi_{j-1}+\phi_j+\phi_{j+1}$. 
   Then $\widetilde{\phi_j}*\phi_j*f=\phi_j*f.$
Using Young's inequality, we have
\begin{equation*}
       {\left\|\F^{-1}_{\re^4}\left[\frac{(\xi,-\eta)}{|\x|+|\eta|}\widehat{\widetilde{\phi}}_j(\x)\widehat{f_j}(\x)\widehat{\widetilde{\phi}}_{j'}(\eta)\widehat{g_{j'}}(\eta)\right]\right\|}_{L^\infty}
      \leq C{\left\|\F^{-1}_{\re^4}\left[\frac{(\xi,-\eta)}{|\x|+|\eta|}\widehat{\widetilde{\phi_j}}(\x)\widehat{\widetilde{\phi_{j'}}}(\eta)\right]\right\|}_{L^1}{\|f_j\|}_{L^\infty}{\|g_{j'}\|}_{L^\infty}.
\end{equation*}
and since $|j-j'|\leq 1$, we have 
\begin{equation*}
   \begin{split}
      {\left\|\F^{-1}_{\re^4}\left[\frac{(\xi,-\eta)}{|\x|+|\eta|}\widehat{\widetilde{\phi}}_j(\x)\widehat{\widetilde{\phi}}_{j'}(\eta)\right]\right\|}_{L^1}
      &={\left\|\F^{-1}_{\re^4}\left[\frac{(\xi,-\eta)}{|\x|+|\eta|}\widehat{\widetilde{\phi}}_{j-j'}(\x)\widehat{\widetilde{\phi_0}}(\eta)\right]\right\|}_{L^1}\\
      &\leq Constant,
   \end{split}
\end{equation*}
where the constant is independent of $j$.
\end{proof}

\begin{lem}
   \label{key}
   Let $T>0$\ and $\theta\in C([0,T];\dot{B}^0_{\infty,1})$\ be a mild solution of (\ref{SQG}).
    Then for any $j\in\Int,\ \pt_t\theta_j\in C([0,T];\dot{B}^0_{\infty,1})$ and for all $t\in(0,T)$, $\theta_j$ satisfies
    $$\pt_t \theta_j + \Lm \theta_j + \N \phi_j * (u\theta) = 0\ \text{in}\ \dot{B}^0_{\infty,1}.$$
\end{lem}

\begin{proof}
   Fix $t\in(0,T)$ and $x\in\re^2$. 
   By definition of the mild solution, for $\phi_j(x-\cdot)\in\S(\re^2)$, we obtain
   $$\phi_j*\th(t)=\phi_j*e^{-t\Lm}\th_0-\N\phi_j*\left(\int_{0}^{t}e^{-(t-s)\Lm}(u\th) ~{\rm d}s\right).$$
   Since $\N\phi_j*$ is a bounded operator in $L^\infty$, we can write
   $$\N\phi_j*\left(\int_{0}^{t}e^{-(t-s)\Lm}(u\th) ~{\rm d}s\right)=\int_{0}^{t}e^{-(t-s)\Lm}(\N\phi_j*(u\th)) ~{\rm d}s,\quad t>0,\quad x\in \re^2.$$
   To prove Lemma \ref{key}, we show that
   $$\pt_t\left(\int_{0}^{t}e^{-(t-s)\Lm}(\N\phi_j*(u\th)) ~{\rm d}s\right)=-\Lm\left(\int_{0}^{t}e^{-(t-s)\Lm}(\N\phi_j*(u\th)) ~{\rm d}s\right)+\N \phi_j * (u\theta).$$
   For $h\in(0,T-t)$,
   \begin{align*}
      &\quad\frac{1}{h}\left(\int_{0}^{t+h}e^{-(t+h-s)\Lm}(\N\phi_j*(u\th)) ~{\rm d}s-\int_{0}^{t}e^{-(t-s)\Lm}(\N\phi_j*(u\th)) ~{\rm d}s\right)\\
      &=\frac{1}{h}\int_{t}^{t+h}e^{-(t+h-s)\Lm}(\N\phi_j*(u\th)) ~{\rm d}s+\frac{1}{h}\int_{0}^{t}(e^{-(t+h-s)\Lm}-e^{-(t-s)\Lm})(\N\phi_j*(u\th)) ~{\rm d}s\\
      &=:I+I\hspace{-1.2pt}I.
   \end{align*}
   The continuity of $e^{-t\Lm}$ and $\N\phi_j*(u\th)$ with respect to time, we have
   $$\lim_{h\to 0}I=\N\phi_j*(u(t)\th(t))\ \text{in}\ \dot{B}^0_{\infty,1}. $$
   By the boundness of $e^{-h\Lm}$ in $\dot{B}^0_{\infty,1}$ and $\Lm\N\phi_j*(u\th)\in C([0,T];\dot{B}^0_{\infty,1})$, we can justify the limit for $I\hspace{-1.2pt}I$ and we obtain
   $$\lim_{h\to 0}I\hspace{-1.2pt}I= \frac{e^{-h\Lm}-Id}{h}\int_{0}^{t}e^{-(t-s)\Lm}(\N\phi_j*(u\th)) ~{\rm d}s\BLACK =-\Lm\left(\int_{0}^{t}e^{-(t-s)\Lm}(\N\phi_j*(u\th)) ~{\rm d}s\right)\ \text{in}\ \dot{B}^0_{\infty,1}.$$
   Therefore $\th_j$ is differentiable in time and
   $$\pt_t \theta_j = - \Lm \theta_j - \N \phi_j * (u\theta)\ \text{in}\ \dot{B}^0_{\infty,1}. $$
\end{proof}

We consider the transport-diffusion equation
\begin{equation}
   \label{TD}
   \begin{cases}
     \pt_t \th + \Lm \th + u \cdot \N \th = f+g,
           &\quad  \   t>0,x\in{\re}^2, \\
     \th(0,x) = \th_0(x),
           &\quad  \ x\in{\re}^2,\\
   \end{cases}
\end{equation}
where $u$\ is a given vector field such that $\div u=0$ and $f,\ g$\ are given external forces.

 The following proposition is a corollary of Theorem 1.2 in \cite{Miao}.
 We give a brief proof. 
\begin{prop}
   \label{M-P}
Let $-1<s<1$. 
Assume that $\th$ is a smooth solution of (\ref{TD}).
Then $\th$ satisfies
\begin{equation}
   {\|\th\|}_{ L^\infty(0,T;\dot{B}^s_{\infty,\infty})}\leq C\exp\big(C{\|\N u\|}_{L^1(0,T;\dot{B}^0_{\infty,\infty})}\big)({\|\th_0\|}_{\dot{B}^s_{\infty,1}}+{\|f\|}_{ L^1(0,T;\dot{B}^s_{\infty,\infty})}+{\|g\|}_{ L^\infty(0,T;\dot{B}^{s-1}_{\infty,\infty})})
\end{equation}
\end{prop}
\begin{proof}
Let $\psi_j$ be the flow of $\displaystyle\sum_{k\leq j-2}u_k$. 
We denote $\overline{f}_j:=f_j\circ\psi_j$. 
Applying $\phi_j*$ to \eqref{TD} and change of variable $x\mapsto \psi_j$, we have 
\begin{equation}
   \label{0429-1}
   \pt_t\overline{\th}_j+\Lm\overline{\th}_j=\overline{f}_j+\overline{g}_j+\overline{R}_j+G_j, 
\end{equation}
where $R_j:=\Big( \displaystyle\sum_{k\leq j-2}u_k-u\Big)\cdot\N\th_j-[\phi_j*,u\cdot\N]\th$ and 
$G_j:=\Lm\overline{\th}_j-(\Lm\th_j)\circ\psi_j$. 
We write $\th_{j,j'}:=\phi_{j'}*\phi_{j}*\th$. 
Applying $\phi_{j'}*$ and using Lemma \ref{0314-1} \eqref{0407-1}, we get 
\begin{equation*}
   \begin{split}
      {\|\overline{\th}_{j,j'}\|}_{L^\infty}&\leq Ce^{-ct2^{j'}}{\|\th_{0,j,j'}\|}_{L^\infty}\\
      &\quad +C\int_{0}^{t}e^{-c(t-\tau)2^{j'}}\big({\|\overline{f}_{j,j'}(\tau)\|}_{L^\infty}+{\|\overline{g}_{j,j'}(\tau)\|}_{L^\infty}+{\|\overline{R}_{j,j'}(\tau)\|}_{L^\infty}+{\|G_{j,j'}(\tau)\|}_{L^\infty}\big) ~{\rm d}\tau.
   \end{split}
\end{equation*}
According to the property of flow and commutator estimate, we obtain 
\begin{equation*}
   \begin{split}
      {\|\overline{\th}_{j,j'}\|}_{L^\infty}&\leq Ce^{-ct2^{j'}}{\|\th_{0,j,j'}\|}_{L^\infty}\\
      &\quad +C2^{j-j'}\int_{0}^{t}e^{-c(t-\tau)2^{j'}}\exp\big(C{\|\N u\|}_{L^1(0,\tau;L^\infty)}\big)\\
      &\hspace{70pt}\times \big({\|f_j(\tau)\|}_{L^\infty}+{\|g_j(\tau)\|}_{L^\infty}+2^{-js}{\|\N u(\tau)\|}_{\dot{B}^0_{\infty,1}}{\|\th(\tau)\|}_{\dot{B}^s_{\infty,\infty}}\\
      &\hspace{90pt}+2^{j'}{\|\N u\|}^\frac{1}{2}_{L^1(0,\tau;L^\infty)}{\|\th_j(\tau)\|}_{L^\infty}\big) ~{\rm d}\tau. 
   \end{split}
\end{equation*}
Taking the $L^\infty$ norm over $[0,t]$ and multiplying both sides by $2^{js}$, we have
\begin{equation*}
   \begin{split}
      &\quad2^{js}{\|\overline{\th}_{j,j'}\|}_{L^\infty}\\
      &\leq C2^{js}{\|\th_{0,j,j'}\|}_{L^\infty}\\
      &\quad+C\max\{1,2^{j-j'}\}2^{j-j'}\exp\big(C{\|\N u\|}_{L^1(0,t;L^\infty)}\big)\big(2^{js}{\|f_j\|}_{L^1(0,t;L^\infty)}+2^{j(s-1)}{\|g_j\|}_{L^\infty(0,t;L^\infty)}\big)\\
      &\quad+C2^{j-j'}\int_{0}^{t}\exp\big(C{\|\N u\|}_{L^1(0,\t;L^\infty)}\big){\|\N u(\tau)\|}_{\dot{B}^0_{\infty,1}}{\|\th(\tau)\|}_{\dot{B}^s_{\infty,\infty}} ~{\rm d}\tau\\
      &\quad+C2^{j-j'}\exp\big(C{\|\N u\|}_{L^1(0,t;L^\infty)}\big){\|\N u\|}^\frac{1}{2}_{L^1(0,t;L^\infty)}2^{js}{\|\th_j\|}_{L^\infty(0,t;L^\infty)}. 
   \end{split}
\end{equation*}
The rest of the proof is the same as \cite{Miao}. 
\end{proof}

Thanks to Lemma \ref{0314-1} \eqref{0407-1}, we obtain convergence of linear solustion.
\begin{lem}
   \label{0407-2}
   Suppose that $\th_0\in \dot{B}^0_{\infty,1}$.
   Then
   $$\lim_{T\to 0}{\|e^{-t\Lm}\th_0\|}_{L^1(0,T;\dot{B}^1_{\infty,1})}=0.$$
\end{lem}

Also next lemma implies convergence of the nonlinear term in the mild solustion of \eqref{SQG}.
\begin{lem}
   \label{0407-3}
   If $\th$ is a mild solustion of \eqref{SQG}, then
   $$\lim_{t\to 0}{\|\th(t)-e^{-t\Lm}\th_0\|}_{\dot{B}^0_{\infty,1}}=0.$$
\end{lem}

Following lemma are valid for functions with high or low frequency restriction.

\begin{lem}
   \label{0314-3}
   Let $s\in\re,\ s'>0$ and $1\leq p,q\leq \infty$.
   Then
   $${\|f_{\text{high}}\|}_{\dot{B}^{s-s'}_{p,q}}\leq C{\|f_{\text{high}}\|}_{\dot{B}^s_{p,q}}\quad \text{and}\quad
   {\|f_{\text{low}}\|}_{\dot{B}^{s+s'}_{p,q}}\leq C{\|f_{\text{low}}\|}_{\dot{B}^s_{p,q}}. $$
\end{lem}

\section{Proof of Theorem \ref{Thm}.}

Let $ w =\theta^{(1)}-\theta^{(2)},\ u^{(i)}= \N ^\perp \Lm^{-1}\th^{(i)}\ (i=1,2)$.
We show that ${\| w\|}_{\dot{B}^0_{\infty,1}}=0$\ in $[0,T_0]$ for some small $T_0>0.$
To that end, we consider the following:
$${\| w_\text{{high}}\|}_{\dot{B}^{-\eta}_{\infty,\infty}}+{\| w_\text{{low}}\|}_{\dot{B}^0_{\infty1}}\ \text{with}\ 0< \eta < 1,$$
where $w_\text{{high}}$ and $w_\text{{low}}$ are defined by \eqref{0321-1}.
\ By Lemma \ref{key},\ for any $j\in\Int,\  w_j=\phi_j * w$\ satisfies
\begin{equation}
   \label{TD2}
   \begin{cases}
      \pt_t  w_j + \Lm  w_j + \N \phi_j * (u^{(1)}\th^{(2)}-u^{(1)}\th^{(2)}) = 0,\\
       w_j(0,x) = 0\\
    \end{cases}
\end{equation}
and
$$u^{(1)}\th^{(1)}-u^{(2)}\th^{(2)}=\frac{1}{2}\sum_{i=1}^{2}((\N ^\perp \Lm^{-1} w)\th^{(i)}+(\N ^\perp \Lm^{-1}\th^{(i)}) w).$$
Moreover, since $\th^{(1)},\th^{(2)}$\ are mild solustions of (\ref{SQG}) and are written by the sum of the linear and the nonlinear parts,
\begin{equation*}
   \begin{split}
     &\frac{1}{2}\sum_{i=1}^{2}((\N ^\perp \Lm^{-1} w)\th^{(i)}+(\N ^\perp \Lm^{-1}\th^{(i)}) w)\\
     &=(\N ^\perp \Lm^{-1} w)e^{-t\Lm}\th_0+(\N ^\perp \Lm^{-1}(e^{-t\Lm}\th_0)) w\\
     &\quad-\frac{1}{2}\sum_{i=1}^{2}((\N ^\perp \Lm^{-1} w)B(u^{(i)}\th^{(i)})+(\N ^\perp \Lm^{-1}B(u^{(i)}\th^{(i)})) w)
   \end{split}
\end{equation*}
where $B$ is defined by \eqref{0313-1}.

First, we consider $ w_\text{{high}}$\ i.e., high frequency part of $ w$.
Applying $\displaystyle\sum_{j'=1}^\infty\phi_{j'}*$\ to (\ref{TD2}), we denote $w_{\text{high},j}:=\phi_j*\displaystyle\sum_{j'=1}^\infty\phi_{j'}*w$ and obtain
\begin{equation}
   \label{0420-1}
   \begin{split}
      &\pt_t  w_{\text{high},j} + \Lm  w_{\text{high},j} + \N \phi_j *  \Big((\N ^\perp \Lm^{-1}(e^{-t\Lm}\th_0)) w \Big)_{\text{high}}\\
      &=-\N \phi_j *\Big((\N ^\perp \Lm^{-1} w)e^{-t\Lm}\th_0-\frac{1}{2}\sum_{i=1}^{2} \big((\N ^\perp \Lm^{-1} w)B(u^{(i)}\th^{(i)})+(\N ^\perp \Lm^{-1}B(u^{(i)}\th^{(i)})) w \big)\Big)_{\text{high}}.
   \end{split}
\end{equation}
Now, we can write
\begin{equation*}
   \begin{split}
       \Big((\N ^\perp \Lm^{-1}(e^{-t\Lm}\th_0)) w \Big)_{\text{high}}&= \Big((\N ^\perp \Lm^{-1}(e^{-t\Lm}\th_0)) w_\text{high} \Big)_{\text{high}}+\Big((\N ^\perp \Lm^{-1}(e^{-t\Lm}\th_0)) w_\text{low}\Big)_{\text{high}}\\
      &=(\N ^\perp \Lm^{-1}(e^{-t\Lm}\th_0)) w_\text{high}\\
      &\quad- \Big((\N ^\perp \Lm^{-1}(e^{-t\Lm}\th_0)) w_\text{high} \Big)_{\text{low}}+ \Big((\N ^\perp \Lm^{-1}(e^{-t\Lm}\th_0)) w_\text{low} \Big)_{\text{high}}.
   \end{split}
\end{equation*}
Thus, we have
\begin{equation}
   \label{0314-10}
   \begin{split}
      &\pt_t  w_{\text{high},j} + \Lm  w_{\text{high},j} + \N \phi_j * ((\N ^\perp \Lm^{-1}(e^{-t\Lm}\th_0)) w_{\text{high}})\\
      &=\N \phi_j * \Bigg( \Big((\N ^\perp \Lm^{-1}(e^{-t\Lm}\th_0)) w_{\text{high}} \Big)_{\text{low}}- \Big((\N ^\perp \Lm^{-1}(e^{-t\Lm}\th_0)) w_{\text{low}} \Big)_{\text{high}}\\
      &\quad-\Big((\N ^\perp \Lm^{-1} w)e^{-t\Lm}\th_0-\frac{1}{2}\sum_{i=1}^{2} \big((\N ^\perp \Lm^{-1} w)B(u^{(i)}\th^{(i)})+(\N ^\perp \Lm^{-1}B(u^{(i)}\th^{(i)})) w \big)\Big)_{\text{high}} \Bigg).
   \end{split}
\end{equation}

Let $0< \eta < 1$.
We divide $w$ into $w_{\text{high}}$ and $w_{\text{low}}$ and consider ${\| w_{\text{high}}\|}_{\dot{B}^{-\eta}_{\infty,\infty}}$\ instead of the $\dot{B}^0_{\infty,1}$\ norm.
Using maximum principle (Proposition \ref{M-P}) to \eqref{0314-10}, $ w_{\text{high}}$ satisfies
\begin{equation}
   \label{high}
   \begin{split}
      &\quad {\| w_{\text{high}}\|}_{ L^\infty(0,T;\dot{B}^{-\eta}_{\infty,\infty})}\\
      &\leq C\exp \Big(C{\|\N e^{-t\Lm}\th_0\|}_{L^1(0,T;\dot{B}^0_{\infty,1})} \Big) \Bigg({ \Big\|\N\cdot \Big((\N ^\perp \Lm^{-1}(e^{-t\Lm}\th_0)) w_{\text{high}} \Big)_{\text{low}} \Big\|}_{ L^1(0,T;\dot{B}^{-\eta}_{\infty,\infty})}\\
      &\quad+{ \Big\|\N\cdot \Big((\N ^\perp \Lm^{-1}(e^{-t\Lm}\th_0)) w_{\text{low}} \Big)_{\text{high}} \Big\|}_{ L^1(0,T;\dot{B}^{-\eta}_{\infty,\infty})}\\
      &\quad+{ \Big\|\N\cdot \Big((\N ^\perp \Lm^{-1} w_{\text{high}})e^{-t\Lm}\th_0 \Big)_{\text{high}} \Big\|}_{ L^1(0,T;\dot{B}^{-\eta}_{\infty,\infty})}\\
      &\quad+{ \Big\|\N\cdot \Big((\N ^\perp \Lm^{-1} w_{\text{low}})e^{-t\Lm}\th_0 \Big)_{\text{high}} \Big\|}_{ L^1(0,T;\dot{B}^{-\eta}_{\infty,\infty})}\\
      &\quad+\sum_{i=1}^{2}{ \Big\|\N\cdot \Big((\N ^\perp \Lm^{-1} w_{\text{high}})B(u^{(i)}\th^{(i)})+(\N ^\perp \Lm^{-1}B(u^{(i)}\th^{(i)})) w_{\text{high}} \Big)_{\text{high}} \Big\|}_{ L^\infty(0,T;\dot{B}^{-\eta-1}_{\infty,\infty})}\\
      &\quad+\sum_{i=1}^{2}{ \Big\|\N\cdot \Big((\N ^\perp \Lm^{-1} w_{\text{low}})B(u^{(i)}\th^{(i)})+(\N ^\perp \Lm^{-1}B(u^{(i)}\th^{(i)})) w_{\text{low}} \Big)_{\text{high}} \Big\|}_{ L^\infty(0,T;\dot{B}^{-\eta-1}_{\infty,\infty})} \Bigg)\\
      &=:C\exp\Big(C{\|\N e^{-t\Lm}\th_0\|}_{L^1(0,T;\dot{B}^0_{\infty,1})}\Big)((I\text{-}1)+(I\text{-}2)+(I\hspace{-1.2pt}I\text{-}1)+(I\hspace{-1.2pt}I\text{-}2)+(I\hspace{-1.2pt}I\hspace{-1.2pt}I\text{-}1)+(I\hspace{-1.2pt}I\hspace{-1.2pt}I\text{-}2)).
   \end{split}
\end{equation}
By Bony's decomposition (\cite{Bony}) , $\N\cdot\N ^\perp \Lm^{-1}w_{\text{high},k}=0$, Lemma \ref{0314-1} \eqref{0314-2} and Lemma~\ref{0314-3}, we write
\begin{equation*}
   \begin{split}
      (I\text{-}1)&\leq C\left({\left\|\left(\sum_{k\leq l-2}(\N ^\perp \Lm^{-1}(e^{-t\Lm}\th_0)_k) w_{\text{high},l}\right)_{\text{low}}\right\|}_{ L^1(0,T;\dot{B}^{-\eta}_{\infty,\infty})}\right.\\
      &\quad + {\left\|\left(\sum_{|k-l|\leq 1}(\N ^\perp \Lm^{-1}(e^{-t\Lm}\th_0)_k) w_{\text{high},l}\right)_{\text{low}}\right\|}_{ L^1(0,T;\dot{B}^{-\eta+1}_{\infty,\infty})}\\
      &\quad +\left. {\left\|\left(\sum_{l\leq k-2}(\N ^\perp \Lm^{-1}(e^{-t\Lm}\th_0)_k) w_{\text{high},l}\right)_{\text{low}}\right\|}_{ L^1(0,T;\dot{B}^{-\eta+1}_{\infty,\infty})}\right),\\
   \end{split}
\end{equation*}

\begin{equation*}
   \begin{split}
      (I\hspace{-1.2pt}I)\text{-}1&\leq C\left({\left\|\left(\sum_{k\leq l-2}(\N ^\perp \Lm^{-1} w_{\text{high},k})(e^{-t\Lm}\th_0)_l\right)_{\text{high}}\right\|}_{ L^1(0,T;\dot{B}^{-\eta+1}_{\infty,\infty})}\right.\\
      &\quad + {\left\|\left(\sum_{|k-l|\leq 1}(\N ^\perp \Lm^{-1} w_{\text{high},k})(e^{-t\Lm}\th_0)_l\right)_{\text{high}}\right\|}_{ L^1(0,T;\dot{B}^{-\eta+1}_{\infty,\infty})}\\
      &\quad + \left.{\left\|\left(\sum_{l\leq k-2}(\N ^\perp \Lm^{-1} w_{\text{high},k})\N(e^{-t\Lm}\th_0)_l\right)_{\text{high}}\right\|}_{ L^1(0,T;\dot{B}^{-\eta}_{\infty,\infty})}\right)\\
   \end{split}
\end{equation*}
and
\begin{equation*}
   \begin{split}
      &\quad(I\hspace{-1.2pt}I\hspace{-1.2pt}I\text{-}1)\\
      &\leq C\sum_{i=1}^{2}\left( {\left\|\left(\sum_{k\leq l-2}\left((\N ^\perp \Lm^{-1} w_{\text{high},k})(B(u^{(i)}\th^{(i)}))_l+(\N ^\perp \Lm^{-1}(B(u^{(i)}\th^{(i)}))_l) w_{\text{high},k}\right)\right)_{\text{high}}\right\|}_{ L^\infty(0,T;\dot{B}^{-\eta}_{\infty,\infty})}\right.\\
      &\quad + {\left\|\N\cdot\left(\sum_{|k-l|\leq 1}\left((\N ^\perp \Lm^{-1} w_{\text{high},k})(B(u^{(i)}\th^{(i)}))_l+(\N ^\perp \Lm^{-1}(B(u^{(i)}\th^{(i)}))_l) w_{\text{high},k}\right)\right)_{\text{high}}\right\|}_{ L^\infty(0,T;\dot{B}^{-\eta-1}_{\infty,\infty})}\\
      &\quad + \left. {\left\|\left(\sum_{l\leq k-2}\left((\N ^\perp \Lm^{-1} w_{\text{high},k})(B(u^{(i)}\th^{(i)}))_l+(\N ^\perp \Lm^{-1}(B(u^{(i)}\th^{(i)}))_l) w_{\text{high},k}\right)\right)_{\text{high}}\right\|}_{ L^\infty(0,T;\dot{B}^{-\eta}_{\infty,\infty})}\right)\\
      &=:(I\hspace{-1.2pt}I\hspace{-1.2pt}I\text{-}1a) +(I\hspace{-1.2pt}I\hspace{-1.2pt}I\text{-}1b) +(I\hspace{-1.2pt}I\hspace{-1.2pt}I\text{-}1c).
   \end{split}
\end{equation*}
Using the product estimate of Besov norm (Lemma \ref{0314-4} and Lemma \ref{0314-5}), we have
\begin{equation}
   \label{0314-6}
   (I\text{-}1)\leq C{\|e^{-t\Lm}\th_0\|}_{L^1(0,T;\dot{B}^0_{\infty,1}\cap\dot{B}^1_{\infty,1})}{\| w_{\text{high}}\|}_{ L^\infty(0,T;\dot{B}^{-\eta}_{\infty,\infty})},
\end{equation}

\begin{equation}
   \label{0315-1}
   (I\hspace{-1.2pt}I\text{-}1)\leq C{\|e^{-t\Lm}\th_0\|}_{L^1(0,T;\dot{B}^1_{\infty,1})}{\| w_{\text{high}}\|}_{ L^\infty(0,T;\dot{B}^{-\eta}_{\infty,\infty})}
\end{equation}
and
$$(I\hspace{-1.2pt}I\hspace{-1.2pt}I\text{-}1a),(I\hspace{-1.2pt}I\hspace{-1.2pt}I\text{-}1c)\leq C\sum_{i=1}^{2}{\|B(u^{(i)}\th^{(i)})\|}_{L^\infty(0,T;\dot{B}^0_{\infty,1})}{\| w_{\text{high}}\|}_{L^\infty(0,T;\dot{B}^{-\eta}_{\infty,\infty})}.$$

In order to  estimate $(I\hspace{-1.2pt}I\hspace{-1.2pt}I\text{-}1b)$, we write
$$(\N ^\perp \Lm^{-1}f)g+(\N ^\perp \Lm^{-1}g)f=\N^\perp((\Lm^{-1}f)g)-(\Lm^{-1}f)(\N^\perp g)+(\N ^\perp \Lm^{-1}g)f. $$
 The Fourier transform in $\re^4$ of the above is written with the separate variables, and\BLACK we have
\begin{equation}
   \label{0314-11}
   \begin{split}
      (\N ^\perp \Lm^{-1}g)f-(\Lm^{-1}f)(\N^\perp g)
      &=\F^{-1}_{\re^4}\left[\eta^\perp\left(\frac{1}{|\eta|}-\frac{1}{|\xi|}\right)\widehat{f}(\xi)\widehat{g}(\eta)\right]\\
      &=\F^{-1}_{\re^4}\left[\frac{(\xi,\eta)\cdot(\xi,-\eta)}{|\xi|+|\eta|}\cdot\frac{1}{|\xi|}\widehat{f}(\xi)\frac{\eta^\perp}{|\eta|}\widehat{g}(\eta)\right]\\
      &=\N_{\re^4}\cdot m(D_1,D_2)((\Lm^{-1}f)(\N ^\perp \Lm^{-1}g)).
   \end{split}
\end{equation}
where $m(D_1,D_2)(fg):=\F^{-1}_{\re^4}\left[\displaystyle \frac{(\xi,-\eta)}{|\xi|+|\eta|}\widehat{f}(\xi)\widehat{g}(\eta)\right]$ and $\N_{\re^4}$ is the gradient \BLACK in $\re^4$.
We replace $f$ and $g$ in \eqref{0314-11} with $ w_{\text{high},k}$ and $(B(u^{(i)}\th^{(i)}))_l$ respectively.
By $\N\cdot\N^\perp((\Lm^{-1}f)g)=0$ and Lemma~\ref{0314-12}, we obtain
\begin{equation}
   \label{0323-2}
   \begin{split}
      &\quad (I\hspace{-1.2pt}I\hspace{-1.2pt}I\text{-}1b)\\
      &\leq C \sum_{i=1}^{2}{\left\|\left(\sum_{|k-l|\leq 1}\N_{\re^4}\cdot m(D_1,D_2)((\Lm^{-1} w_{\text{high},k})(\N ^\perp \Lm^{-1}(B(u^{(i)}\th^{(i)}))_l))\right)_{\text{high}}\right\|}_{ L^\infty(0,T;\dot{B}^{-\eta}_{\infty,\infty})}\\
      &\leq C\sum_{i=1}^{2}\sup_{j\leq 0}2^{(-\eta+1)j}\sum_{k\geq j-4}\sum_{|k-l|\leq 1}{\|\Lm^{-1} w_{\text{high},k}\|}_{L^\infty(0,T;L^\infty)}{\|(B(u^{(i)}\th^{(i)}))_l\|}_{L^\infty(0,T;L^\infty)}\\
      &\leq C\sum_{i=1}^{2}{\|B(u^{(i)}\th^{(i)})\|}_{L^\infty(0,T;\dot{B}^0_{\infty,1})}{\| w_{\text{high}}\|}_{L^\infty(0,T;\dot{B}^{-\eta}_{\infty,\infty})}
   \end{split}
\end{equation}
and
\begin{equation}
   \label{0315-2}
   (I\hspace{-1.2pt}I\hspace{-1.2pt}I\text{-}1)\leq C\sum_{i=1}^{2}{\|B(u^{(i)}\th^{(i)})\|}_{L^\infty(0,T;\dot{B}^0_{\infty,1})}{\| w_{\text{high}}\|}_{L^\infty(0,T;\dot{B}^{-\eta}_{\infty,\infty})}.
\end{equation}

By Lemma \ref{0314-1} \eqref{0314-2}, Lemma \ref{0323-1} and Lemma \ref{0314-3}, we estimate $(I\text{-}2)$ as follows:
\begin{equation}
   \label{0314-7}
   \begin{split}
      (I\text{-}2)&\leq C{ \Big\|\Big((\N ^\perp \Lm^{-1}(e^{-t\Lm}\th_0)) w_{\text{low}} \Big)_{\text{high}} \Big\|}_{ L^1(0,T;\dot{B}^1_{\infty,\infty})}\\
      &\leq C{\|e^{-t\Lm}\th_0\|}_{L^1(0,T;\dot{B}^0_{\infty,1}\cap\dot{B}^1_{\infty,1})}{\| w_{\text{low}}\|}_{L^\infty(0,T;\dot{B}^0_{\infty,1}\cap\dot{B}^1_{\infty,1})}\\
      &\leq C {\|e^{-t\Lm}\th_0\|}_{L^1(0,T;\dot{B}^0_{\infty,1}\cap\dot{B}^1_{\infty,1})}{\| w_{\text{low}}\|}_{L^\infty(0,T;\dot{B}^0_{\infty,1})}.
   \end{split}
\end{equation}
Similarly, we obtain
\begin{equation}
   \label{0314-8}
   (I\hspace{-1.2pt}I\text{-}2)\leq C{\|e^{-t\Lm}\th_0\|}_{L^1(0,T;\dot{B}^0_{\infty,1}\cap\dot{B}^1_{\infty,1})}{\| w_{\text{low}}\|}_{L^\infty(0,T;\dot{B}^0_{\infty,1})}.
\end{equation}
Using Lemma \ref{0314-3} and a similar argument as the estimate $(I\hspace{-1.2pt}I\hspace{-1.2pt}I\text{-}1b)$ \eqref{0323-2} , we obtain
\begin{equation}
   \label{0314-9}
   \begin{split}
      (I\hspace{-1.2pt}I\hspace{-1.2pt}I\text{-}2)&\leq C\sum_{i=1}^{2}{ \Big\|\N\cdot \Big((\N ^\perp \Lm^{-1} w_{\text{low}})B(u^{(i)}\th^{(i)})+(\N ^\perp \Lm^{-1}B(u^{(i)}\th^{(i)})) w_{\text{low}} \Big)_{\text{high}} \Big\|}_{ L^\infty(0,T;\dot{B}^{-1}_{\infty,\infty})}\\
      &\leq C\sum_{i=1}^{2}{\|B(u^{(i)}\th^{(i)})\|}_{L^\infty(0,T;\dot{B}^0_{\infty,1})}{\| w_{\text{low}}\|}_{L^\infty(0,T;\dot{B}^0_{\infty,1})}.
   \end{split}
\end{equation}
From \eqref{0314-6}, \eqref{0315-1} and \eqref{0315-2}-\eqref{0314-9}, we obtain
\begin{equation}
   \label{high2}
   \begin{split}
      &\quad{\| w_{\text{high}}\|}_{ L^\infty(0,T;\dot{B}^{-\eta}_{\infty,\infty})}\\
      &\leq C\exp(C{\|\N e^{-t\Lm}\th_0\|}_{L^1(0,T;\dot{B}^0_{\infty,1})}) \Bigg({\|e^{-t\Lm}\th_0\|}_{L^1(0,T;\dot{B}^0_{\infty,1}\cap\dot{B}^1_{\infty,1})}+\sum_{i=1}^{2}{\|B(u^{(i)}\th^{(i)})\|}_{L^\infty(0,T;\dot{B}^0_{\infty,1})} \Bigg)\\
      &\quad \times({\| w_{\text{high}}\|}_{ L^\infty(0,T;\dot{B}^{-\eta}_{\infty,\infty})}+{\| w_{\text{low}}\|}_{ L^\infty(0,T;\dot{B}^0_{\infty,1})}).
   \end{split}
\end{equation}

Next, we consider $ w_\text{{low}}$. 
A similar argument as high frequency estimate (\ref{0420-1}) and using the Duhamel's principle, we have
\begin{equation}
   \begin{split}
      \begin{split}
         &\pt_t  w_{\text{low},j} + \Lm  w_{\text{low},j} + \N \phi_j *  \Big((\N ^\perp \Lm^{-1}(e^{-t\Lm}\th_0)) w \Big)_{\text{low}}\\
         &=-\N \phi_j *\Big((\N ^\perp \Lm^{-1} w)e^{-t\Lm}\th_0-\frac{1}{2}\sum_{i=1}^{2} \big((\N ^\perp \Lm^{-1} w)B(u^{(i)}\th^{(i)})+(\N ^\perp \Lm^{-1}B(u^{(i)}\th^{(i)})) w \big)\Big)_{\text{low}}.
      \end{split}
   \end{split}
\end{equation}
and 
\begin{equation}
   \begin{split}
      \begin{split}
         w_{\text{low},j}(t) = - \int_{0}^{t} e^{-(t-s)\Lm} \Bigg(\N \phi_j *  &\Big((\N ^\perp \Lm^{-1}(e^{-t\Lm}\th_0)) w + (\N ^\perp \Lm^{-1} w)e^{-t\Lm}\th_0\\
         & + \frac{1}{2}\sum_{i=1}^{2} \big((\N ^\perp \Lm^{-1} w)B(u^{(i)}\th^{(i)})+(\N ^\perp \Lm^{-1}B(u^{(i)}\th^{(i)})) w \big)\Big)_{\text{low}}\Bigg) ~{\rm d}s.
      \end{split}
   \end{split}
\end{equation}
By Lemma \ref{0314-1} \eqref{0407-1}, $e^{-(t-s)\Lm}$ is a bounded operator in $\dot{B}^0_{\infty,1}$ and we obtain
\begin{equation}
   \begin{split}
      &\quad{\| w_{\text{low}}\|}_{ L^\infty(0,T;\dot{B}^0_{\infty,1})}\\
      &\leq C\Bigg({\Big\|\N\cdot\Big((\N ^\perp \Lm^{-1}(e^{-t\Lm}\th_0)) w_{\text{high}}+(\N ^\perp \Lm^{-1} w_{\text{high}})e^{-t\Lm}\th_0\Big)_{\text{low}}\Big\|}_{ L^1(0,T;\dot{B}^0_{\infty,1})}\\
      &\quad+{\Big\|\N\cdot\Big((\N ^\perp \Lm^{-1}(e^{-t\Lm}\th_0)) w_{\text{low}}+(\N ^\perp \Lm^{-1} w_{\text{low}})e^{-t\Lm}\th_0\Big)_{\text{low}}\Big\|}_{ L^1(0,T;\dot{B}^0_{\infty,1})}\\
      &\quad+\sum_{i=1}^{2}{\Big\|\N\cdot\Big((\N ^\perp \Lm^{-1} w_{\text{high}})B(u^{(i)}\th^{(i)})+(\N ^\perp \Lm^{-1}B(u^{(i)}\th^{(i)})) w_{\text{high}}\Big)_{\text{low}}\Big\|}_{ L^1(0,T;\dot{B}^0_{\infty,1})}\\
      &\quad+\sum_{i=1}^{2}{\Big\|\N\cdot\Big((\N ^\perp \Lm^{-1} w_{\text{low}})B(u^{(i)}\th^{(i)})+(\N ^\perp \Lm^{-1}B(u^{(i)}\th^{(i)})) w_{\text{low}}\Big)_{\text{low}}\Big\|}_{ L^1(0,T;\dot{B}^0_{\infty,1})} \Bigg).\\
   \end{split}
\end{equation}
Since derivatives for low frequency component can be estimated by a constant, we also have
\begin{equation}
   \label{low2}
   \begin{split}
      &\quad{\| w_{\text{low}}\|}_{ L^\infty(0,T;\dot{B}^0_{\infty,1})}\\
      &\leq C\Bigg({\|e^{-t\Lm}\th_0\|}_{L^1(0,T;\dot{B}^0_{\infty,1}\cap\dot{B}^1_{\infty,1})}+\sum_{i=1}^{2}{\|B(u^{(i)}\th^{(i)})\|}_{L^1(0,T;\dot{B}^0_{\infty,1})} \Bigg)\\
      &\quad \times({\| w_{\text{high}}\|}_{ L^\infty(0,T;\dot{B}^{-\eta}_{\infty,\infty})}+{\| w_{\text{low}}\|}_{ L^\infty(0,T;\dot{B}^0_{\infty,1})}).
   \end{split}
\end{equation}

Summing up (\ref{high2}) and (\ref{low2}), since Lemma \ref{0407-2} and Lemma \ref{0407-3}, for any $\epsilon>0$, there exist $0<T_0 \ll 1$ such that
$${\| w_{\text{high}}\|}_{ L^\infty(0,T_0;\dot{B}^{-\eta}_{\infty,\infty})}+{\| w_{\text{low}}\|}_{ L^\infty(0,T_0;\dot{B}^0_{\infty,1})}\leq C\epsilon({\| w_{\text{high}}\|}_{ L^\infty(0,T_0;\dot{B}^{-\eta}_{\infty,\infty})}+{\| w_{\text{low}}\|}_{ L^\infty(0,T_0;\dot{B}^0_{\infty,1})}).$$
By taking $\epsilon$ sufficiently small, we obtain
$${\| w_{\text{high}}\|}_{ L^\infty(0,T_0;\dot{B}^{-\eta}_{\infty,\infty})}+{\| w_{\text{low}}\|}_{ L^\infty(0,T_0;\dot{B}^0_{\infty,1})}\leq 0.$$
We conclude that
$$\th^{(1)}=\th^{(2)}\ \text{in}\ \dot{B}^0_{\infty,1}\ \text{for}\ t\in[0,T_0].$$

We prove that ${\| w\|}_{\dot{B}^0_{\infty,1}}=0$\ in the entire interval $[0,T]$, by contradiction argument.
Set
$$\t^*:=\sup\{\t\in[0,T)\ |\ {\|\th^{(1)}(t,\cdot)-\th^{(2)}(t,\cdot)\|}_{\dot{B}^0_{\infty,1}}=0\ \text{for}\ t\in[0,\t]\}.$$
Assume that $\t^*<T$. 
From the continuity in time of $\th^{(1)}$ and $\th^{(2)}$, we have $\th^{(1)}(\t^*)=\th^{(2)}(\t^*).$
By the same argument as local uniqueness , there exists $\del'>0$ such that $\th^{(1)}(s+\t^*)=\th^{(2)}(s+\t^*)$ for $0\leq s\leq\del'$.
This is a contradiction with the definition of $\t^*$.
Therefore $\t^*=T$.
\begin{flushright}
   $\square$
\end{flushright}

\noindent
{\bf Data availability statement.} This manuscript has no associated data.\\
{\bf Conflict of Interest.} The author declares that he has no conflict of interest.

\end{document}